\documentclass[reqno,final]{amsart}
\usepackage{natbib}  
\usepackage{fancyhdr} 
\usepackage{color} 
\usepackage{hyperref} 
\usepackage{graphicx} 

\usepackage{pstricks}
\usepackage{amssymb}
\usepackage{amsmath}
\usepackage{graphicx}
\usepackage{multirow}
\usepackage{float}
\usepackage{dsfont}

\definecolor{aleacolor}{rgb}{0.16,0.59,0.78}

\hypersetup{
breaklinks,
colorlinks=true,
linkcolor=aleacolor,
urlcolor=aleacolor,
citecolor=aleacolor}

\renewcommand{\headheight}{11pt}
\renewcommand{\cite}{\citet}

\theoremstyle{plain}
\newtheorem{theorem}{Theorem}[section]                                          
\newtheorem{proposition}[theorem]{Proposition}                          
\newtheorem{lemma}[theorem]{Lemma}
\newtheorem{corollary}[theorem]{Corollary}

\theoremstyle{definition}

\theoremstyle{remark}

\makeatletter \@addtoreset{equation}{section} \makeatother
\renewcommand\theequation{\thesection.\arabic{equation}}
\renewcommand\thefigure{\thesection.\arabic{figure}}
\renewcommand\thetable{\thesection.\arabic{table}}

\begin{document}

\title[Stopping Times of Random Walks on a Hypercube]{Stopping Times of Random Walks on a Hypercube}

\author{Cl\'audia Peixoto}
\author{Diego Marcondes}

\address{Instituto de Matem\'atica e Estat\'istica \newline
Universidade de S\~ao Paulo, Brazil}

\email{claudiap@ime.usp.br, dmarcondes@ime.usp.br}
\urladdr{\url{http://www.ime.usp.br/~dmarcondes}}


\begin{abstract}
  A random walk on a $N$-dimensional hypercube is a discrete time stochastic process whose state space is the set $\{-1,+1\}^{N}$, which has uniform probability of reaching any neighbour state, and probability zero of reaching a non-neighbour state, in one step. This random walk is often studied as a process associated with the Ehrenfest Urn Model. This paper aims to present results about the time that such random walk takes to self-intersect and to return to a set of states. We also present results about the time that the random walk on a hypercube takes to visit a given set and a random set of states. Asymptotic distributions and bounds are presented for these times. The coupling of random walks is widely used as a tool to prove the results.
\end{abstract}

\maketitle

\section{Introduction}

A random walk on a $N$-dimensional hypercube is a discrete time stochastic process whose state space is the $N$-dimensional hypercube, i.e., the set $\{-1,+1\}^{N}$. There are $2^{N}$ possible states for this random walk, which may be seem as vertices of the hypercube. The random walk on the hypercube is given by flipping one coordinate (or spin) at each time and is often studied as a process associated with the Ehrenfest Urn Model (see \cite{ehr}), when there is an interest not only on the proportion of particles at each urn, but also on where each particle is at each time. Indeed, we may consider the ball $i$ to be in the first (second) urn when the configuration of the \textit{i-th} vertex of the hypercube assumes value $1$ ($-1$).

Although extensively studied and applied, the random walk on a hypercube still lacks some simple and important results, especially about its stopping times. Therefore, this paper aims to present results about stopping times of random walks on a $N$-dimensional hypercube. We treat the first self-intersection time, the time to revisit the path taken from $0$ to $[N^{\gamma}], 0 < \gamma < 1$, the time to visit a given set and the time to visit a random set.

The first theorem shows that the first self-intersection of a random walk on a $N$-dimensional hypercube is, with probability 1, as $N \rightarrow \infty$, a two-step one, i.e., a return characterized by two consecutive spins flips on the same vertex. From this result, it immediately follows that the distribution of the time of such self-intersection, properly standardized, converges to an exponential distribution with rate $1$.

Propositions \ref{prop1} through \ref{prop4} and Theorem \ref{T3} treat the return to the first $[N^{\gamma}], 0 < \gamma < 1,$ states visited by the random walk, establishing bounds and asymptotic distributions for the time of such return. They also present results for the time that it takes to visit a fixed set $V \subset \{-1,+1\}^{N}$. Finally, Proposition \ref{prop5} and Theorem \ref{T4} present results for the time the random walk takes to first visit a random set $M_{N} \subset \{-1,+1\}^{N}$ with a properly chosen average size.

For simplicity, Theorem \ref{T3} and Propositions \ref{prop1} through \ref{prop4} are enunciated for the aperiodic random walk, while Theorems \ref{T2} and \ref{T4}, and Proposition \ref{prop5}, are enunciated for the periodic random walk, although they can be immediately extended to the periodic and aperiodic one, respectively. Indeed, the aperiodic process takes the same path of the periodic one, that is, the sequence of states is the same for both processes, but the aperiodic one is slower, as it may remain at a state for a while. Therefore, based on this consideration, the results of this paper may in fact be extended to both processes. 

Asymptotic results for random walks on a hypercube have been widely studied. The time that a particle takes to reach its stationary distribution was showed by \cite{asy} to occur around $\frac{1}{4} N \log N$, and the total variation distance to stationariness at this threshold was studied. The probability of hitting a vertex $a$ before hitting a vertex $b$, whenever $a$ and $b$ shared the same edge, starting at any position, was presented by \cite{volkov}. The structure of the set of unvisited sets was studied by \cite{cooper} and \cite{mat}; results about its transition probabilities were given by \cite{letac} and \cite{sco}; and some interesting applications of random walks on a hypercube were presented by \cite{crowe} and \cite{gilbert}. Hitting times of the Ehrenfest chain were determined by \cite{mat2} by the use of coupled random walks on a hypercube. In \cite{nestoridi2017} the author study a non-local random walk on the hypercube, which at each step flips $k$ randomly chosen coordinates, and her main result states the mixing time of this random walk. However, the results of this paper about the stopping times of random walks on the hypercube seems to not have been treated before in this literature, so our results come as a complement to the vast theory about random walks on a hypercube. 

This paper contains the results of the master's thesis \cite{tese}, and extends some results of \cite{tese1}, in which the hitting time of the random set $M_{N}$ is studied in the presence of a ``\textit{temperature}'' in each vertex which characterizes the waiting time of the random walk at the vertex. In Section \ref{sec2} we present the notation used throughout the paper. In Section \ref{sec3} we present the main results of the paper, and in Section \ref{sec4} we present their proofs.

\section{Random Walk on a Hypercube}
\label{sec2}

A $N$-dimensional hypercube is denoted by $H_{N}=\{-1,+1\}^{N}$ and its vertices are represented by $\eta$, i.e., $\eta=\{\eta_{1},\dots,\eta_{N}\}, \eta_{i} \in \{-1,+1\}, i \in \{1, \dots, N\}$. Given $j \in \{1, \dots ,N\}$ and $\eta \in H_{N}$, the vertex $\eta^{j}$ obtained from $\eta$ by a spin flip at its $j$th coordinate is given by

\begin{equation}
\label{1}
(\eta^{j})_{i}= 
\begin{cases}
\eta_{i}, & \text{if } i \neq j;\\
-(\eta_{i}), & \text{if } i = j;
\end{cases}, \ i \in \{1, \dots, N\}.
\end{equation}

A random walk on a hypercube is a stochastic process whose state space is the hypercube $H_{N}$ and whose transition from one state to another is made by a spin at a vertex of the hypercube, as defined by (\ref{1}). The transition probabilities of this process may be defined in two ways, dividing these random walks in two types: aperiodic and periodic.

The aperiodic random walk is denoted by $\sigma$, where $\sigma(t) \in H_{N}$ is its state at a time $t \in \mathbb{N}$. The random walk $\sigma$ may be obtained from two sequences of independent random variables $I(t)$ and $U(t), \ t \in \mathbb{N}$, defined in a probability space $(\Omega_{N},\mathbb{F}_{N},\mathbb{P}_{N})$. The random variables $I(t)$ assume values in $\{1, \dots, N\}$, are independent, identically distributed and, for any $k \in \{1, \dots, N\}$, $\mathbb{P}_{N}\{I(t)=k\}=1/N$, and the random variables $U(t)$ are independent and identically distributed, with uniform distribution in $[0,1]$. Then, for all $t \in \mathbb{N}$, $\omega \in \Omega_{N}$ and $i \in \{1, \dots, N\}, \ \sigma_{i}(t)$ is given by	
\[
\sigma_{i}(t,\omega)=
\begin{cases}
\sigma_{i}(t-1,\omega), & \text{if } I(t,\omega) \neq i;\\
+1, & \text{if } I(t,\omega)=i; \ U(t,\omega) < \frac{1}{2};\\
-1, & \text{if } I(t,\omega)=i; \ U(t,\omega) \geq \frac{1}{2}.
\end{cases}
\]

On the other hand, the periodic random walk is denoted by $\xi$, where $\xi(t) \in H_{N}$ is its state at a time $t \in \mathbb{N}$. The random walk $\xi$, defined on the probability space $(\Omega_{0},\mathbb{F}_{0},\mathbb{P}_{0})$, is a random walk with transition probabilities given by	
\[
\mathbb{P}_{0}\Big(\xi(k+1)=\eta^{i}\Big|\xi(k)=\eta\Big)=\frac{1}{N}, \forall k \geq 0; \eta \in H_{N},
\]
where $\eta^{i}$ is one of the $N$ states that may be reached from $\eta$ by a spin, as defined in (\ref{1}), $i \in \{1, \dots, N\}$. Note that $\xi$ may be also defined as a function of the sequence of random variables $I(t), t \in \mathbb{N}$. The difference between $\xi$ and $\sigma$ is that $\xi$ changes its state at each step with probability 1, while $\sigma$ has probability $1/2$ of not changing its state. Both random walks are used to prove our results, as some proofs are more elegant for $\xi$ and others for $\sigma$.

The random walks with initial state $(-1, \dots, -1)$ are denoted by $\xi^{-}$ and $\sigma^{-}$, those with initial state $(+1, \dots, +1)$ are denoted by $\xi^{+}$ and $\sigma^{+}$ and the random walks with initial state $\eta \in H_{N}$ are denoted by $\xi^{\eta}$ and $\sigma^{\eta}$. When the initial state has no importance, the random walks are denoted simply by $\xi$ and $\sigma$. The indexes of the probability spaces defined above are omitted if there is no doubt about which one is being referred.

Coupled random walks $\sigma^{\eta}$ and $\sigma^{\varsigma}$ are constructed using the same $\omega$, i.e., the same index  $I(t,\omega)$ and the same value of $U(t,\omega)$, given $\sigma^{\eta}(t),\sigma^{\varsigma}(t)$ and $I(t+1)=i$, in the following way:
\[
\begin{cases}
\text{if} \ U(t+1) <\frac{1}{2}, & \text{ then } \sigma^{\eta}_{i}(t+1)=+1,\  \sigma^{\varsigma}_{i}(t+1)=+1;\\
\text{if} \ U(t+1) \geq \frac{1}{2}, & \text{ then } \sigma^{\eta}_{i}(t+1)=-1,\  \sigma^{\varsigma}_{i}(t+1)=-1.
\end{cases}
\]

The distance at time $t$ between the two random walks $\sigma^{\eta}$ e $\sigma^{\varsigma}$ is defined by
\[
D_{N}^{\eta, \varsigma}(t)=\frac{1}{2N}\sum_{i=1}^{N}|\sigma_{i}^{\eta}(t)-\sigma_{i}^{\varsigma}(t)|.
\] 
Note that $D_{N}^{\eta, \varsigma}(t)$ is an Ehrenfest Model in $\bigg\{ 0,\frac{1}{N},...,D_{N}^{\eta, \varsigma}(0)-\frac{1}{N},D_{N}^{\eta, \varsigma}(0) \bigg\}$, in which zero is an absorbing state.

The time taken by coupled aperiodic random walks $\sigma^{+}$ and $\sigma^{-}$ to meet is defined by $t^{-}_{N}=\min(t>0:\sigma^{+}(t)=\sigma^{-}(t))$. A well-known result for $t^{-}_{N}$ is that
\begin{equation}
\lim\limits_{N\rightarrow \infty}\frac{t^{-}_{N}}{N\log N} = 1 \text{ in probability. }
\label{res1}
\end{equation}
See \cite[Proposition III.1]{tese1} for a proof of this result.

Let $\sigma^{\eta}$ and $\sigma^{\varsigma}$ be coupled random walks on $H_{N}$. Suppose that $D^{\eta,\varsigma}_{N}(0)=[Nf]/N, 0 \leq f \leq 1$. Then, it follows from (\ref{res1}) that
\begin{equation}
\label{res3}
\lim\limits_{N \rightarrow \infty} \mathbb{P}\Big(\sigma^{\eta}(t(N))\neq \sigma^{\varsigma}(t(N))\Big)=0,
\end{equation}
for any $t(N)$ that satisfies $\lim\limits_{N \rightarrow \infty}\frac{t(N)}{NlogN}=\infty$.

The result below, that also follows from (\ref{res1}), presents an upper bound for the rate of convergence to the equilibrium:
\begin{equation}
\label{res4}
\lim\limits_{N \rightarrow \infty}\bigg| \mathbb{P}(\sigma^{+}(t(N))=\eta) - \frac{1}{2^{N}} \bigg| = 0,
\end{equation}
for all $\eta \in H_{N}$ and $t(N)$ satisfying $\lim\limits_{N \rightarrow \infty}\frac{t(N)}{NlogN}=\infty$.

\section{Results}
\label{sec3}

Our first result treats the time of the first self-intersection of $\xi$, that is given by the random variable $S_{N} =\min(t \geq 2:\xi(t) \in \{\xi(0), \dots, \xi(t-1)\})$, and the time of a $2l$-step return of $\xi$, i.e, the first return to a state in $2l$ steps, $l \geq 1$. Note that such return happens if, and only if, a $2k$-step return ($k < l$) did not happen to such state (or any other) and a return happened at the $2l$-th step. A $2l$-step first return may be defined as a function of the random variables $\{I(1), \dots, I(2l)\}$ and occurs if, and only if,
\begin{enumerate}
	\item In all the vectors $\Big\{(I(m), \dots, I(m+s)): m = 1, \dots, 2l; m+s \leq 2l\Big\}$ there is one value $j \in \{1, \dots, N\}$ which appears an odd number of times. This guarantees that $\xi(t) \neq \xi(m)$ for all $m < t \leq m + s$, so that no $2k$-step return has occurred for $k < l$.
	\item If $j \in \{I(1), \dots, I(2l)\}$ then $j$ appears an even number of times in it. This guarantees that a $2l$-step return has occurred, as all the vertices are the same as the initial state.
\end{enumerate}

In order to establish if a $2l$-step return has occurred, we may apply convenient functions to a sample of $\{I(1),\dots,I(2l)\}$. To this purpose, let $\mathbf{x}=(i_{1}, \dots, i_{2l}) \in \{1, \dots, N\}^{2l}$, $\mathbf{x}^{[j,k]}=(i_{j}, \dots, i_{k}), 1 \leq j < k \leq 2l$ and $\mathds{1}_{(A)}(B)$ be the usual Kronecker's delta. We define for all $\mathbf{x} \in \{1, \dots, N\}^{2l}$ the functions $f_{2l}$, $h_{2l}$ and $g_{2l}$ as
\[
\begin{cases}
f_{2l}(\mathbf{x})= & \prod\limits_{j=1}^{N}\bigg[\mathds{1}_{\{0,2, \dots, 2l\}}\bigg(\sum\limits_{k=1}^{2l}\mathds{1}_{\{j\}}(i_{k})\bigg) \bigg]; \\
\\
h_{2l}(\mathbf{x})= & \prod\limits_{j=1}^{2l}\Bigg[\prod\limits_{\substack{k=1 \\ (j,k) \neq (1,l)}}^{\lfloor \frac{2l+1-j}{2} \rfloor}\bigg(1 - f_{2k}\Big(\mathbf{x}^{[j,j+2k-1]}\Big)  \bigg)  \Bigg], l>1; \\
\\
g_{2l}(\mathbf{x})= &

\begin{cases}
f_{2l}(\mathbf{x}) h_{2l}(\mathbf{x}), & \text{if } l>1; \\
f_{2l}(\mathbf{x}), & \text{if } l=1;
\end{cases}
.
\end{cases}
\]

Note that if functions $h_{2l}(\textbf{x})$ and $f_{2l}(\textbf{x})$, in which $\textbf{x}$ is a sample of $\{I(1), \dots,\\ I(2l)\}$, are equal to 1, then conditions 1 and 2, respectively, are being satisfied by the sample, so that if $g_{2l}(\textbf{x}) = 1$ then a return in $2l$-steps happened at the sample. Therefore, defining the set $J_{l} \subset \{1, \dots, N\}^{2l}$ as
\[
J_{l}= \big\{(i_{1}, \dots, i_{2l}) \in \{1, \dots, N\}^{2l}: g_{2l}((i_{1}, \dots, i_{2l}))=1 \big\},
\]
we have that the time of the first $2l$-step return of the random walk $\xi$ may be defined by
\[
\Gamma_{l} =  \min\big(t \geq 2l: \big(I(t-(2l-1)), \dots, I(t)\big) \in J_{l}\big).
\]
Our first result states that, as $N \rightarrow \infty$, the first return of $\xi$ is a 2-step return with probability 1.

\begin{theorem}
	\label{T2}
	\[
	\lim\limits_{N \rightarrow \infty} \mathbb{P}(S_{N}=\Gamma_{1})=1.
	\]	
\end{theorem}

From this theorem it follows that, as $N \rightarrow \infty$, the number of distinct states already visited by a random walk $\xi$ at time $[N^{\gamma}$], $0 < \gamma < 1$, is $[N^{\gamma}]+1$, with probability 1.
\begin{corollary}
	For $0 < \gamma < 1$,
	\[
	\lim\limits_{N \rightarrow \infty} \mathbb{P}\Big(\Big|\{\xi(0), \dots, \xi([N^{\gamma}])\}\Big| = [N^{\gamma}] + 1\Big)=1.
	\]
	\label{col1}
\end{corollary}

Still from Theorem \ref{T2}, it follows that $\frac{S_{N}}{N}$ converges weakly to an exponential law.

\begin{corollary}
	The random variable $N^{-1}S_{N}$ converges weakly to a mean one exponential law.
	\label{col2}
\end{corollary}

We now treat the return of $\sigma^{+}$ to its first $[N^{\gamma}], 0 < \gamma < 1$, visited states, i.e, the path taken from $0$ to $[N^{\gamma}]$, as denoted by $V(0,[N^{\gamma}])=\{\sigma^{+}(0), \dots, \sigma^{+}([N^{\gamma}])\}$. The first return of $\sigma^{+}$ to $V(0,[N^{\gamma}])$ and the first visit of $\sigma^{\eta}$ to $V(0,[N^{\gamma}])$ may be defined, respectively, by $R_{N} = \min (t > N^{\gamma}: \sigma^{+}(t) \in V(0,[N^{\gamma}]))$ and $R^{\eta}_{N} = \min (t>0: \sigma^{\eta}(t) \in V(0,[N^{\gamma}]))$. Note that $\sigma^{+}$ and $\sigma^{\eta}$ are coupled. Consider from now on that $N^{\gamma}=[{N^{\gamma}}]$ and denote $V(0,[N^{\gamma}])$ by $V$. Propositions \ref{prop1} through \ref{prop4} give bounds to $R_{N}$ and $R^{\eta}_{N}$.

First, we note that, as $N \rightarrow \infty$, $N^{1+\delta}$, $0 < \delta < 1/2$, is a lower bound for $R_{N}$ and $R_{N}^{\eta}$.
\begin{proposition}
	\label{prop1}
	For $0 < \gamma < 1$, $0 < \delta < 1/2$ and $\forall \eta \notin V$,
	\begin{equation*}
	\lim\limits_{N \rightarrow \infty} \mathbb{P}(R_{N} > N^{1+\delta}) = \lim\limits_{N \rightarrow \infty} \mathbb{P}(R^{\eta}_{N} > N^{1+\delta})= 1
	\end{equation*}
	
\end{proposition}

The lower bound above may be improved to any bound $t(N)$ satisfying $\lim\limits_{N \rightarrow \infty}\frac{t(N)N^{\gamma}}{2^{N}}=0$.

\begin{proposition}
	\label{prop3}
	For $0 < \gamma < 1$ and for all $t(N)$ satisfying $\lim\limits_{N \rightarrow \infty}\frac{t(N)N^{\gamma}}{2^{N}}=0$,
	\[
	\lim\limits_{N \rightarrow \infty} \mathbb{P}(R_{N} > t(N)) = 1.
	\]
\end{proposition}

We note that any bound $t(N)$ satisfying $\lim\limits_{N \rightarrow \infty}\frac{t(N)^{1-\epsilon}\nu(V)}{N(logN + 1)}=0$, in which $\epsilon > 0$ and $\nu(\cdot)$ is the uniform measure in $H_{N}$, i.e., $\nu(A) = \frac{|A|}{2^{N}}, \forall A \subset H_{N}$, is an upper bound for $R_{N}$.

\begin{proposition}
	\label{prop4}
	For $\epsilon > 0$ and for all $t(N)$ satisfying $\lim\limits_{N \rightarrow \infty}\frac{t(N)^{1-\epsilon}\nu(V)}{N(logN + 1)}=\infty$, in which $\nu(\cdot)$ is the uniform measure in $H_{N}$,
	\[
	\lim\limits_{N \rightarrow \infty} \mathbb{P}(R_{N} \leq t(N)) = 1.
	\]
\end{proposition}

We now engage in determining the weak convergence of $R_{N}$. For this purpose, define $\beta_{N} = \min (t \in \mathbb{N}: \mathbb{P}(R_{N} \geq t) \leq e^{-1})$. The next theorem states that the distribution of $R_{N}$, standardized by $\beta_{N}$, converges to a mean one exponential law.

\begin{theorem}
	\label{T3}
	For $0 < \gamma < 1$ and $t > 0$,
	\[
	\lim\limits_{N \rightarrow \infty} \mathbb{P} \Bigg(\frac{R_{N}}{\beta_{N}} > t \Bigg) = e^{-t}.
	\]
\end{theorem}

Finally, we engage in finding the hitting time for a random set $ M_{N} \subset H_{N}$, defined on the probability space $(\bar{\Omega},\bar{\mathbb{F}},\bar{\mathbb{P}})$. Each vertex of the hypercube is in $M_{N}$ with probability $1/N^{\gamma}, \gamma > 0$, independently of each other. Let the time that the process $\xi$ takes to reach the set $M_{N}$ be defined as $\Theta = \min (t > 0: \xi(t) \in M_{N})$. Proposition \ref{prop5} states that the expected value of the probability of $\Theta$ being greater than $N^{\gamma}t$ equals the survival function of a mean one exponential law.

\begin{proposition}
	\label{prop5}For $0 < \gamma < 1$ and $t >0$,
	\[
	\lim\limits_{N \rightarrow \infty} \bar{E}(\mathbb{P}(\Theta > N^{\gamma}t)) = e^{-t}.
	\]
\end{proposition}

Lastly, we note that, not only the expected value of the probability of $\Theta$ being greater than $N^{\gamma}t$ equals the survival function of a mean one exponential law, but also the limiting distribution of $\frac{\Theta}{N^{\gamma}}$ is such exponential law.

\begin{theorem}
	\label{T4}
	For $0 < \gamma < 1$ and $\epsilon, t > 0$,
	\[
	\lim\limits_{N \rightarrow \infty} \bar{\mathbb{P}}\Big(|\mathbb{P}(\Theta> N^{\gamma}t) - e^{-t}| > \epsilon \Big) = 0.
	\]
\end{theorem}

\section{Proofs}
\label{sec4}

A proof for the theorems, propositions and corollaries above are presented in this section. Some lemmas are enunciated and proved in order to assist the proofs of Theorems \ref{T2} and \ref{T3}.\\

\begin{lemma}
	\label{l1}
	For $l \geq 3$
	\[
	\mathbb{P}\Big( \Gamma_{l} \leq n\Big) \leq \frac{8n}{N^{3}}.
	\]
\end{lemma}

\begin{proof}
	First, we show that $\mathbb{P}\Big( \{I(1), \dots, I(2l)\} \in J_{l}\Big) \leq \frac{8}{N^{3}}$. We have that
	\[
	\mathbb{P}\Big( \{I(1), \dots, I(2l)\} \in J_{l}\Big) = \frac{|J_{l}|}{N^{2l}}.
	\]
	Now, note that $|J_{l}| \leq 2N^{2l-2}$, because, given the first $2l-2$ coordinates, the last two coordinates of a vector $v \in J_{l}$ are fixed, but a permutation. Furthermore, the set $J_{l}$ may be divided into two disjoint sets: its vectors in which the last three coordinates are distinct from one another and those in which only two of the last three coordinates are distinct. Those sets are denoted by $J_{l}'$ and $J_{l}''$, respectively. Thus, $|J_{l}|=|J_{l}'| + |J_{l}''|$ and	
	
	\textit{a)} $|J_{l}'| \leq 3! \ N^{2l-3}$, because, given the first $2l-3$ coordinates, the last three coordinates must be fixed, but a permutation, if the vector is in $J_{l}'$.
	
	\textit{b)} $|J_{l}''| \leq N \ |J_{l-1}|$, because if the pair that is at the last three coordinates is disregarded, exactly one return of the type $|J_{l-1}|$ happens; furthermore, the pair may assume $N$ values.
	
	Therefore,	
	\begin{align*}
	\dfrac{|J_{l}|}{N^{2l}} \leq \frac{3! \ N^{2l-3} + N \ |J_{l-1}|}{N^{2l}} \leq \frac{3! \ N^{2l-3} + N \ (2N^{2l-4})}{N^{2l}} = \frac{8}{N^{3}}.
	\end{align*}
	and it follows that	
	\begin{align*}
	& \mathbb{P} \Big(\Gamma_{l} \leq n \Big)  = \sum\limits_{k=2l}^{n} \mathbb{P} \Big(\Gamma_{l} = k \Big) \\ 
	& = \sum\limits_{k=2l}^{n} \! \mathbb{P} \! \Big(\ \! \! (I(j\!-\!2l\!+1), \! \dots, \! I(j)) \! \notin \! J_{l}, \! j < k; \! (I(k-2l+1), \! \dots, \! I(k)) \! \in \! J_{l} \! \Big) \\
	& \leq (n - 2l) \mathbb{P}\Big((I(1), \dots, I(2l)) \in J_{l} \Big) \leq \frac{8n}{N^{3}}.
	\end{align*}	
\end{proof}

\begin{proof}[Proof of Theorem \ref{T2}]
	First of all, note that $S_{N}= \min\limits_{l \geq 1} \ \Gamma_{l}$, and it is enough to show that	
	\[
	\lim\limits_{N \rightarrow \infty} \mathbb{P} \Bigg(\Gamma_{1} < N^{1+\delta} < \min\limits_{2 \leq l \leq \frac{N^{1+\delta}}{2}} \Gamma_{l} \Bigg) = 1
	\]
	for some $\delta > 0$. On the one hand,	
	\begin{align*}
	\lim\limits_{N \rightarrow \infty} \mathbb{P} \Big(\Gamma_{1} < N^{1+\delta} \Big) & = \lim\limits_{N \rightarrow \infty} 1 - \mathbb{P} \Big(\Gamma_{1} \geq N^{1+\delta} \Big) \\ & = \lim\limits_{N \rightarrow \infty} 1 - \Big(1 - \frac{1}{N} \Big)^{N^{1 + \delta}} = 1.
	\end{align*}
	
	On the other hand, applying Lemma \ref{l1},
	\begin{align}
	\nonumber
	\mathbb{P} \Bigg(\min\limits_{2 \leq l \leq \frac{N^{1+\delta}}{2}} \Gamma_{l} > N^{1+\delta} \Bigg) & = 1 - \mathbb{P} \Bigg(\min\limits_{2 \leq l \leq \frac{N^{1+\delta}}{2}} \Gamma_{l} \leq N^{1+\delta} \Bigg) \\ \label{eq20}
	& \geq 1 - N^{1+\delta} \Bigg(\frac{2}{N^{2}}\Bigg) - \sum\limits_{l=3}^{\frac{N^{1+\delta}}{2}} \frac{8N^{1+\delta}}{N^{3}}
	\end{align}
	and for $0 < \delta < \frac{1}{2}$ the limit of (\ref{eq20}) is 1 as $N \rightarrow \infty$. 		
\end{proof}

\begin{proof}[Proof of Corollary \ref{col1}]
	\begin{align*}
	\lim\limits_{N \rightarrow \infty} \mathbb{P} \Big(\Big|\{\xi(0), \dots, \xi(N^{\gamma})\} \Big| = N^{\gamma} + 1 \Big) & = \lim\limits_{N \rightarrow \infty} \mathbb{P} \Big(S_{N}> N^{\gamma} \Big) \\ & = \lim\limits_{N \rightarrow \infty} \mathbb{P} \Big(\ \Gamma_{1} > N^{\gamma} \Big) \\
	& = \lim\limits_{N \rightarrow \infty} \Bigg(1-\frac{1}{N} \Bigg)^{N^{\gamma}} = 1.  
	\end{align*}	
\end{proof}

\begin{proof}[Proof of Corollary \ref{col2}]
	It is enough to show that, for some $t>0$,
	\[
	\lim\limits_{N \rightarrow \infty} \mathbb{P} \Big(S_{N} > tN \Big)=e^{-t}.
	\]
	
	Note that
	\begin{align*}	
	\Big| \mathbb{P} \Big(S_{N} > tN \Big) - \mathbb{P} \Big(\Gamma_{1} & > tN \Big)   \Big| \leq \mathbb{P} \Big(\Gamma_{1} \neq S_{N} \Big).
	\end{align*}
	
	From Theorem \ref{T2}, $\lim\limits_{N \rightarrow \infty} \mathbb{P} \Big(\Gamma_{1} \neq S_{N} \Big) = 0$, and	
	\[
	\lim\limits_{N \rightarrow \infty} \mathbb{P} \Big( \Gamma_{1} > tN \Big) = \lim\limits_{N \rightarrow \infty} \Big(1 - \frac{1}{N} \Big)^{Nt}=e^{-t}
	\]
	so the result follows.
\end{proof}

\begin{proof}[Proof of Proposition \ref{prop1}]
	From Theorem \ref{T2},
	\[
	\lim\limits_{N \rightarrow \infty} \mathbb{P}\Big( \cup_{l \geq 2} \Big[ \Gamma_{l} < N^{1 + \delta} \Big] \Big) = 0.
	\]
	
	Therefore,
	\begin{align*}
	\lim\limits_{N \rightarrow \infty} \mathbb{P}(R_{N} \leq N^{1+\delta}) = \lim\limits_{N \rightarrow \infty} \mathbb{P} ( \sigma^{+}(N^{\gamma} + 1) = \sigma^{+}(N^{\gamma})) \leq \lim\limits_{N \rightarrow \infty} \frac{N^{\gamma}}{N} = 0
	\end{align*}
	so that it is enough to show that
	\[
	\lim\limits_{N \rightarrow \infty} \Big|\mathbb{P}(R_{N} > N^{1 + \delta}) - \mathbb{P}(R_{N}^{\eta} > N^{1 + \delta})\Big| = 0.
	\]
	
	However, it follows from (\ref{res3}) that	
	\begin{align*}
	& \Big|\mathbb{P}(R_{N} > N^{1 + \delta}) - \mathbb{P}(R_{N}^{\eta} > N^{1 + \delta})\Big| \leq\\
	& \leq \mathbb{P}(R_{N} > N^{1 + \delta}, R_{N}^{\eta} \leq N^{1 + \delta}) + \mathbb{P}(R_{N}^{\eta} > N^{1 + \delta}, R_{N} \leq N^{1 + \delta}) \\
	& \leq \sup\limits_{\eta_{0} \in H_{N}} \mathbb{P}(\sigma^{\eta_{0}}(N^{1 + \delta}) \neq \sigma^{+}(N^{1 + \delta})) \xrightarrow[]{N \rightarrow \infty} 0.
	\end{align*}	
\end{proof}

\begin{proof}[Proof of Proposition \ref{prop3}]
	We have by Proposition \ref{prop1} that
	\begin{align*}
	& \mathbb{P}(R_{N}  \leq t(N)) \leq \mathbb{P}(R_{N} \leq N^{1+\delta}) + \mathbb{P}(N^{1+\delta} < R_{N} \leq t(N)) \\
	\leq & \mathbb{P}(R_{N} \leq N^{1+\delta}) + \mathbb{P}(\sigma^{+}(t) \in V, \text{ for some } N^{1+\delta} < t \leq t(N)) \\
	\leq & o(N) + \sum\limits_{\eta \in H_{N}} \nu(\eta) |\mathbb{P}(\sigma^{+}(t) \in V, \text{ for some } N^{1+\delta} < t \leq t(N)) - \\
	& - \mathbb{P}(\sigma^{\eta}(t) \in V, \text{ for some } N^{1+\delta} < t \leq t(N))| + \\ 
	& + \sum\limits_{\eta \in H_{N}} \nu(\eta) \mathbb{P}(\sigma^{\eta}(t) \in V, \text{ for some } N^{1+\delta} < t \leq t(N)) \\
	\leq & o(N) +\sum\limits_{\eta \in H_{N}} \nu(\eta) \sup\limits_{\eta_{0} \in H_{N}} \mathbb{P}(\sigma^{+}(N^{1+\delta}) \neq \sigma^{\eta_{0}}(N^{1+ \delta})) + \\ & + \sum\limits_{\eta \in H_{N}} \nu(\eta)  \sum\limits_{u=N^{1+\delta}}^{t(N)} \mathbb{P}(\sigma^{\eta}(u) \in V) \\
	\leq & o(N) \! + \! \! \sum\limits_{\eta \in H_{N}} \! \nu(\eta) \! \sup\limits_{\eta_{0} \in H_{N}} \! \mathbb{P}(\sigma^{+}(N^{1+\delta}) \! \neq \! \sigma^{\eta_{0}}(N^{1+ \delta})) \! + \! \frac{(N^{\gamma} \! + \! 1)t(N)}{2^{N}}
	\end{align*}
	The limit above, as $N \rightarrow \infty$, is zero by the hypothesis $\lim\limits_{N \rightarrow \infty} \frac{t(N)N^{\gamma}}{2^{N}} = 0$ and (\ref{res3}).	
\end{proof}

\begin{proof}[Proof of Proposition \ref{prop4}]
	Consider the events $A_{s}=\{\sigma^{\eta}(s) \in V\}$ and the random variable $Z= \sum\limits_{s=0}^{t(N)} \mathds{1}_{\{A_{s}\}}$. Note that $\{Z > 0\} = \{R_{N}^{\eta} \leq t(N)\}$. Applying Paley-Zygmund inequality, for $\delta > 0$,	
	\begin{align*}
	& \mathbb{P}(Z > 0) \geq \frac{\Big[E\Big(\sum\limits_{s=N^{1+ \delta}}^{t(N)} \mathds{1}_{\{A_{s}\}}\Big)\Big]^{2}}{E\Big(\Big[\sum\limits_{s=N^{1+ \delta}}^{t(N)} \mathds{1}_{\{A_{s}\}}\Big]^{2}\Big)} + o(N)\\
	& = \frac{\Big[\sum\limits_{s=N^{1+ \delta}}^{t(N)} \mathbb{P}(A_{s})\Big]^{2}}{E\Big(\sum\limits_{s=N^{1+ \delta}}^{t(N)}(\mathds{1}_{\{A_{s}\}})^{2} + \sum\limits_{u \neq s}(\mathds{1}_{\{A_{u}\}})(\mathds{1}_{\{A_{s}\}})\Big)} + o(N) \\
	& = \frac{\Bigg[\sum\limits_{s=N^{1+ \delta}}^{t(N)} \frac{|V|}{2^{N}}\Bigg]^{2}}{\sum\limits_{s=N^{1+ \delta}}^{t(N)}\frac{|V|}{2^{N}} + \sum\limits_{u \neq s} \mathbb{P}(A_{u} \cap A_{s})} + o(N) \\
	& = \frac{(t(N)-N^{1+ \delta})^{2}(\nu(V))^{2}}{(t(N)-N^{1+ \delta})\nu(V) \!+ \! \! \! \! \! \! \! \! \! \! \sum\limits_{|u-s|>N^{1 + \delta}} \! \! \! \! \!  \mathbb{P}(A_{u} \cap A_{s}) + \! \! \! \! \! \! \! \! \! \! \sum\limits_{0 < |u-s|\leq N^{1 + \delta}} \! \! \! \! \! \mathbb{P}(A_{u} \cap A_{s})} + o(N).
	\end{align*}
	
	On the one hand,
	\begin{align*}
	& \sum\limits_{|u-s|>N^{1 + \delta}} \mathbb{P}(A_{u} \cap A_{s})  = \sum\limits_{k=N^{1+ \delta}}^{t(N)} (t(N)-k+1) \mathbb{P}(\sigma(k) \in V, \sigma(0) \in V) \\
	\leq & \sum\limits_{k=N^{1+ \delta}}^{t(N)} (t(N)-k+1) \sum\limits_{\eta \in V} \nu(\eta) \mathbb{P}(\sigma(k) \in V | \sigma(0) = \eta) \\ 
	= & \sum\limits_{k=N^{1+ \delta}}^{t(N)} (t(N)-k+1) \times \Bigg\{ \sum\limits_{\eta \in V} \nu(\eta) \sum\limits_{\xi \in H_{N}} \nu(\xi) \Big[\mathbb{P} (\sigma(k) \in V | \sigma(0) = \eta) - \\ 
	& - \mathbb{P}(\sigma(k) \in V | \sigma(0)= \xi)\Big] +  \sum\limits_{\eta \in V} \nu(\eta) \sum\limits_{\xi \in H_{N}} \nu(\xi) \mathbb{P}(\sigma(k) \in V | \sigma(0) = \xi) \Bigg\} \\	
	\leq & \sum\limits_{k=N^{1+ \delta}}^{t(N)} (t(N)-k+1) \Bigg[\sum\limits_{\eta \in V} \nu(\eta) \sum\limits_{\xi \in H_{N}} \nu(\xi) \mathbb{P}(\sigma^{\eta}(k) \neq \sigma^{\xi}(k)) + \nu(V)^{2} \Bigg] \\
	\leq & \sum\limits_{k=N^{1+ \delta}}^{t(N)} (t(N)-k+1) \Big[ \nu(V) \frac{N(\log N+1)}{k} + \nu(V)^{2} \Big] \\
	\leq & \sum\limits_{k= 1}^{t(N)} \frac{t(N)}{k} \nu(V) N(\log N + 1) + \sum\limits_{k=1}^{t(N)} t(N)\nu(V)^{2} \\
	\leq & t(N) \log(t(N))  N (\log(N) + 1) \nu(V) + t(N)^{2} \nu(V)^{2} \\
	\leq & t(N)^{1+ \epsilon} N (\log N + 1) \nu(V) + t(N)^{2} \nu(V)^{2}. \\
	\end{align*}
	
	On the other hand,	
	\begin{align*}
	\sum\limits_{0 < |u-s| \leq N^{1+ \delta}} \mathbb{P}(A_{u} \cap A_{s}) & = \sum\limits_{k=1}^{N^{1+\delta}} (t(N)-k+1) \mathbb{P}(\sigma^{\eta}(k) \in V, \sigma(0) \in V) \\
	& \leq \sum\limits_{k=1}^{N^{1+\delta}} t(N) \nu(V) = N^{1+\delta} t(N) \nu(V). \\
	\end{align*}
	
	Therefore,	
	\begin{align*}
	\lim\limits_{N \rightarrow \infty} \mathbb{P}(Z > 0) = \lim\limits_{N \rightarrow \infty} \mathbb{P}(R^{\eta}_{N} \leq t(N)) = 1, \\
	\end{align*}
	by hypothesis. Now, it is enough to prove that
	\[
	\lim\limits_{N \rightarrow \infty} |\mathbb{P}(R_{N}^{\eta} > t(N)) - \mathbb{P}(R_{N} > t(N))| = 0.
	\]
	
	By Proposition \ref{prop1} and (\ref{res3}),	
	\begin{align}
	\nonumber
	|\mathbb{P}(R_{N}^{\eta} > t(N)) - \mathbb{P}(R_{N} > t(N))|& \leq \mathbb{P}(R_{N}^{\eta} > t(N), R_{N} \leq t(N)) + \\ \nonumber & + \mathbb{P}(R_{N} > t(N), R_{N}^{\eta} \leq t(N)) \\ \nonumber
	& \leq \mathbb{P}(R_{N}^{\eta} > N^{1+ \delta}, R_{N} \leq N^{1+ \delta}) + \\ \nonumber & + \mathbb{P}(R_{N} > N^{1+ \delta}, R_{N}^{\eta} \leq N^{1+ \delta}) \\ 
	& \leq \sup\limits_{\eta_{0} \in H_{N}} \mathbb{P}(\sigma^{\eta_{0}}(N^{1+ \delta}) \neq \sigma^{+}(N^{1+ \delta})) \label{rp4}
	\end{align}
	and as $N \rightarrow \infty$, (\ref{rp4}) goes to zero.		
\end{proof}

\begin{lemma}
	\label{l3}
	\[
	\lim\limits_{N \rightarrow \infty} \mathbb{P}(R_{N} \geq \beta_{N}) = e^{-1}.
	\]
\end{lemma}

\begin{proof}
	By the definition of $\beta_{N}$,
	\[
	\mathbb{P}(R_{N} \geq \beta_{N}) \leq e^{-1} < \mathbb{P}(R_{N} \geq \beta_{N} - 1).
	\]
	
	As
	\[
	0 \leq \mathbb{P}(R_{N} \geq \beta_{N} - 1) - \mathbb{P}(R_{N} \geq \beta_{N}) \leq \mathbb{P}(\beta_{N} - 1 \leq R_{N} < \beta_{N}),
	\]
	the proof may be completed applying the Markov property. Note that
	
	\begin{align*}
	&\mathbb{P}(\beta_{N} - 1 \leq R_{N} < \beta_{N}) = \\ & =  \mathbb{P}(\beta_{N} - 1 \leq R_{N} < \beta_{N} | \sigma^{+}(\beta_{N} - 1) \notin V) \times \mathbb{P}(\sigma^{+}(\beta_{N} - 1) \notin V) \\
	& =  \mathbb{P}( \sigma(1) \in V | \sigma(0) \notin V) \times \Big[\mathbb{P}(\sigma^{+}(\beta_{N} - 1) \notin V) - \\ & - \sum\limits_{\eta \in H_{N}} \nu(\eta) \mathbb{P}(\sigma^{\eta}(\beta_{N}-1) \notin V) + \sum\limits_{\eta \in H_{N}} \nu(\eta) \mathbb{P}(\sigma^{\eta}(\beta_{N}-1) \notin V)\Big] \\
	& \leq \sum\limits_{\eta \in H_{N}} \nu(\eta)\Big[\mathbb{P}(\sigma^{+}(\beta_{N} - 1) \notin V) - \mathbb{P}(\sigma^{\eta}(\beta_{N}-1) \notin V)\Big] + \Bigg(\frac{N^{\gamma}(2^{N} - |V|)}{N2^{N}}\Bigg) \\
	& \leq \sum\limits_{\eta \in H_{N}} \nu(\eta) \sup\limits_{\eta_{0} \in H_{N}} \mathbb{P}(\sigma^{+}(\beta_{N}-1) \neq \sigma^{\eta_{0}}(\beta_{N} - 1)) + \Bigg(\frac{N^{\gamma}(2^{N} - |V|)}{N2^{N}}\Bigg).	
	\end{align*}
	
	As $N \rightarrow \infty$, the first term is zero by (\ref{res3}) and the second is zero by Corollary \ref{col1} and the hypothesis $0 < \gamma < 1$. Therefore,
	\[
	\lim\limits_{N \rightarrow \infty}\mathbb{P}(R_{N} \geq \beta_{N}) = e^{-1}.
	\]		
\end{proof}

\begin{lemma}
	\label{l4}
	For each integer $n > 0$, there exists an $\alpha$ satisfying $e^{-1} \leq \alpha < 1$ such that, for every $N > M_{n} > 0$,
	\[
	\mathbb{P}(R_{N} \geq n\beta_{N}) \leq \alpha^{n}.
	\]
\end{lemma}

\begin{proof}
	The Lemma is proved by induction. For $n = 1$ the result is immediate by definition. Now assume that the inequality holds good for the integer $n$. Then, applying the Markov property, for a $N$ big enough, i.e., $N > M_{n} > 0$,
	\begin{align*}
	\mathbb{P}(R_{N} \geq \beta_{N}(n + 1)) & = \sum\limits_{\eta \notin V} \mathbb{P}(R_{N} \geq \beta_{N}n, \sigma(\beta_{N}n) = \eta) \times \mathbb{P}(R_{N} \geq \beta_{N} | \sigma(0) = \eta) \\
	& \leq \mathbb{P}(R_{N} \geq \beta_{N}n) \sup\limits_{\eta_{0} \notin V} \mathbb{P}(R_{N} \geq \beta_{N} | \sigma(0) = \eta_{0}) \\
	& \leq \alpha^{n} \sup\limits_{\eta_{0} \notin V} \mathbb{P}(R_{N} \geq \beta_{N} | \sigma(0) = \eta_{0}) \\
	& \leq \alpha^{n} e^{-1} \leq \alpha^{n+1}.
	\end{align*}	
\end{proof}

\begin{proof}[Proof of Theorem \ref{T3}]
	To show that $R_{N}$ standardized by $\beta_{N}$ converges to an exponential law with rate $1$ as $N \rightarrow \infty$, it is enough to prove that
	\begin{enumerate}
		\item[\textit{(1)}] $\lim\limits_{N \rightarrow \infty} \Big| \mathbb{P}(R_{N} > \beta_{N}(t + s)) - \mathbb{P}(R_{N}> \beta_{N}t)\mathbb{P}(R_{N} > \beta_{N}s) \Big| = 0.$
		\item[\textit{(2)}] $\lim\limits_{N \rightarrow \infty} E\Big(\frac{R_{N}}{\beta_{N}}\Big) = 1.$
	\end{enumerate} 
	
	Note that the first item guarantees that if the law $R_{N}/\beta_{N}$ converges when $N \rightarrow \infty$, then this limit must be an exponential law (maybe a degenerate one). On the other hand, Lemma \ref{l3} along with this item implies that if $t$ is a positive rational number, then the limit
	\[
	\lim\limits_{N \rightarrow \infty} \mathbb{P}(R_{N} \geq \beta_{N}t)
	\]
	exists and equals $e^{-t}$. As the exponential law is continuous, it is enough to show the convergence for all $t \in \mathbb{R}$, what concludes the proof. The second item guarantees that such exponential law has rate 1. This technique was applied in \cite{tese2}, where more details are presented.
	
	\textit{Proof of (1):} First, the result is proved for an initial state chosen uniformly from $H_{N}$. For this end, note the following three facts:
	\begin{enumerate}
		\item[(a)] 
		\begin{align*}
		& \Big|\sum\limits_{\eta \in H_{N}} \nu(\eta) \mathbb{P}(R_{N}^{\eta} > \beta_{N}(t + s)) - \\ & -\!\sum\limits_{\eta \in H_{N}} \! \nu(\eta) \! \mathbb{P}(\sigma^{\eta}(u) \! \notin \! V, \forall \! u \in \{1, \dots, \beta_{N}t\} \cup \{\beta_{N}t+N^{1 + \delta}, \dots, \beta_{N}(t + s)\}) \Big| \\
		& \leq \sum\limits_{\eta \in H_{N}} \nu(\eta) \mathbb{P}(\sigma^{\eta}(u) \in V, \text{ for some } u \in \{\beta_{N}t+1, \dots, \beta_{N}t + N^{1 + \delta}\}) \\
		& \leq \sum\limits_{\eta \in H_{N}} \nu(\eta) \sum\limits_{u=\beta_{N}t + 1}^{\beta_{N}t + N^{1 + \delta}} \sum\limits_{\varsigma \in V} \mathbb{P}(\sigma^{\eta}(u) = \varsigma) \leq  \sum\limits_{\eta \in H_{N}} \nu(\eta) \sum\limits_{u=\beta_{N}t}^{\beta_{N}t + N^{1 + \delta}} \sum\limits_{\varsigma \in V} \frac{1}{2^{N}} \\
		& = \frac{N^{1 + \delta}(N^{\gamma} + 1)}{2^{N}} \xrightarrow[]{N \rightarrow \infty} 0.
		\end{align*}
		
		\item[(b)]
		\begin{align*}
		& \Big|\sum\limits_{\eta \in H_{N}} \nu(\eta) \mathbb{P}(R^{\eta}_{N} > \beta_{N}s) - \sum\limits_{\eta \in H_{N}} \nu(\eta) \mathbb{P}(\sigma^{\eta}(u) \notin V, \forall u \in \{N^{1 + \delta}, \dots, \beta_{N}s\})\Big| \\
		& \leq  \sum\limits_{\eta \in H_{N}} \nu(\eta) \mathbb{P}(\sigma^{\eta}(u) \in V, \text{ for some } u \in \{1, \dots, N^{1 + \delta}\}) \\
		& \leq  \sum\limits_{\eta \in H_{N}} \nu(\eta) \sum\limits_{u=1}^{N^{1 + \delta}} \sum\limits_{\varsigma \in V} \mathbb{P}( \sigma^{\eta}(u) = \varsigma) \leq \frac{N^{1 + \delta}(N^{\gamma}+1)}{2^{N}} \xrightarrow[]{N \rightarrow \infty} 0.
		\end{align*}
		
		\item[(c)]
		Applying the Markov property, facts (a) and (b), and (\ref{res3}), we have that
		\begin{align*}
		& \Big|\sum\limits_{\eta \in H_{N}} \nu(\eta) \mathbb{P}\Big(R^{\eta}_{N} > \beta_{N}(t+s)\Big) - \sum\limits_{\eta \in H_{N}} \nu(\eta) \mathbb{P}\Big(R^{\eta}_{N} > \beta_{N}t\Big)  \sum\limits_{\eta \in H_{N}} \nu(\eta) \mathbb{P}\Big(R^{\eta}_{N} > \beta_{N}s\Big)\Big| \\
		\leq & \Big|\sum\limits_{\eta \in H_{N}} \sum\limits_{\kappa \notin V} \nu(\eta) \mathbb{P}\Big(R_{N}^{\eta} > \beta_{N}t, \sigma^{\eta}(\beta_{N}t) = \kappa\Big) \\
		& \times \Big[\mathbb{P}\Big(\sigma^{\kappa}(u) \notin V, \forall u \in \{N^{1 + \delta}, \dots, \beta_{N}s\}\Big) - \\ & - \mathbb{P}\Big(\sigma^{\eta}(u) \notin V, \forall u \in \{N^{1 + \delta}, \dots, \beta_{N}s\}\Big)  \Big]\Big|  \\ \leq & \sum\limits_{\eta \in H_{N}} \nu(\eta) \sum\limits_{\kappa \notin V} \sup\limits_{\kappa_{0} \in H_{N}} \mathbb{P} \Big(\sigma^{\kappa_{0}}(N^{1 + \delta}) \neq \sigma^{\eta}(N^{1 + \delta})\Big) \xrightarrow[]{N \rightarrow \infty} 0.
		\end{align*}
	\end{enumerate}
	
	Now, it is enough to show that
	\[
	\lim\limits_{N \rightarrow \infty}	\Big|\mathbb{P}(R_{N} > \beta_{N}t) - \mathbb{P}(R_{N}^{\eta} > \beta_{N}t) \Big| = 0.
	\]
	
	Note that, by Proposition \ref{prop1},
	\begin{align*}
	& \Big|\mathbb{P}(R_{N} > \beta_{N}t) - \mathbb{P}(\sigma^{+}(u) \notin V, \forall u \in \{N^{1+\delta}, \dots, \beta_{N}t\}) \Big| \leq \\ \leq & \Big|\mathbb{P}(R_{N} \leq \beta_{N}t) - \mathbb{P}(N^{1 + \delta} \leq R_{N} \leq \beta_{N}t) \Big| \\
	= & \mathbb{P}(R_{N} < N^{1 + \delta}) \xrightarrow[]{N \rightarrow \infty} 0
	\end{align*}
	and, analogously,	
	\begin{align*}
	\Big|\mathbb{P}(R_{N}^{\eta} > \beta_{N}t)  - \mathbb{P}(\sigma^{\eta}(u) \notin V, \forall u \in \{N^{1+\delta}, \dots, \beta_{N}t\}) \Big|
	\leq \mathbb{P}(R_{N}^{\eta} < N^{1 + \delta}) \xrightarrow[]{N \rightarrow \infty} 0.
	\end{align*}
	
	Therefore,	
	\begin{align}
	\nonumber
	& \Big|\mathbb{P}(R_{N} > \beta_{N}t) - \mathbb{P}(R_{N}^{\eta} > \beta_{N}t)\Big| \\ \nonumber
	\leq & \Big|\mathbb{P}(\sigma^{+}(u) \notin V, \forall u \in \{N^{1+\delta}, \dots, \beta_{N}t\}) - \mathbb{P}(\sigma^{\eta}(u) \notin V, \forall u \in \{N^{1+\delta}, \dots, \beta_{N}t\}) \Big| \\ \label{23}
	\leq & \sup\limits_{\eta_{0} \in H_{N}} \mathbb{P}(\sigma^{+}(N^{1 + \delta}) \neq \sigma^{\eta_{0}}(N^{1 + \delta})).	
	\end{align}
	
	Applying (\ref{res3}), as $N \rightarrow \infty$, (\ref{23}) is zero and we have the result.
	
	\textit{Proof of (2):}
	By definition,
	\begin{align*}
	E\Bigg(\frac{R_{N}}{\beta_{N}}\Bigg) = \int_{0}^{\infty} \mathbb{P}\Bigg(\frac{R_{N}}{\beta_{N}} > t \Bigg) dt.
	\end{align*}
	
	As $N \rightarrow \infty$, because of Lemma \ref{l4}, the Lebesgue's Dominated Convergence Theorem may be applied and by Lemma \ref{l3}		
	\[
	\lim\limits_{N \rightarrow \infty} E\Bigg(\frac{R_{N}}{\beta_{N}}\Bigg) = \int_{0}^{\infty} \lim\limits_{N \rightarrow \infty} \mathbb{P}\Bigg(\frac{R_{N}}{\beta_{N}} > t \Bigg) dt = 1.
	\]
\end{proof}

\begin{proof}[Proof of Proposition \ref{prop5}]
	First, note that, for a given $\bar{\omega} \in \bar{\Omega}$,
	\[
	\mathbb{P}(\Theta > N^{\gamma}t) = \mathds{1}_{\{\eta \notin M_{N}\}} \frac{1}{N} \sum\limits_{i_{1} = 1}^{N} \Bigg(\mathds{1}_{\{\eta^{i_{1}} \notin M_{N}\}} \dots \frac{1}{N} \sum\limits_{i_{N^{\gamma}t} = 1}^{N} \mathds{1}_{\{\eta^{i_{1}, \dots, i_{N^{\gamma}t}} \notin M_{N}\}} \Bigg).
	\]
	
	Let
	\[
	F_{1}(N) \!= \! \Big\{ \!(i_{1}, \! \! \dots, \! \! i_{N^{\gamma}t}) \! \in \! \{1, \! \dots, \! N\}^{N^{\gamma}t}; l = 1, \! \dots, \! N^{\gamma}t: \! \eta^{i_{1}, \! \dots, \! i_{l}} \! \notin \! \{\eta, \! \eta^{i_{1}}, \! \dots, \! \eta^{i_{1}, \! \dots, \! i_{l-1}}\} \! \Big\}
	\]
	and $F_{2}(N) = \{1, \dots, N\}^{N^{\gamma}t}/F_{1}(N)$.
	
	If $i'=(i_{1}, \dots, i_{N^{\gamma}t})$, then	 
	\begin{align*}
	\bar{E}\Big(\mathbb{P}(\Theta > N^{\gamma}t)\Big) = & \frac{1}{N^{N^{\gamma}t}} \sum\limits_{i' \in F_{1}(N)} \bar{E}\Big(\mathds{1}_{\{\eta^{i_{1}} \notin M_{N}\}} \dots \mathds{1}_{\{\eta^{i_{1}, \dots, i_{N^{\gamma}t}} \notin M_{N}\}}\Big) \\
	& + \frac{1}{N^{N^{\gamma}t}} \sum\limits_{i' \in F_{2}(N)} \bar{E}\Big(\mathds{1}_{\{\eta^{i_{1}} \notin M_{N}\}} \dots \mathds{1}_{\{\eta^{i_{1}, \dots, i_{N^{\gamma}t}} \notin M_{N}\}}\Big) \\
	= & \frac{|F_{1}(N)|}{N^{N^{\gamma}t}} \Bigg(1 - \frac{1}{N^{\gamma}}\Bigg)^{N^{\gamma}t+1} + \frac{|F_{2}(N)|}{N^{N^{\gamma}t}} \Bigg(1 - \frac{1}{N^{\gamma}}\Bigg)^{N^{\gamma}t+1}.
	\end{align*}
	
	Applying Corollary \ref{col1},
	\[
	\lim\limits_{N \rightarrow \infty } \frac{|F_{1}(N)|}{N^{N^{\gamma}t}} = 1 \text{ and } \lim\limits_{N \rightarrow \infty } \frac{|F_{2}(N)|}{N^{N^{\gamma}t}} = 0.
	\]
	
	Therefore,
	\[
	\lim\limits_{N \rightarrow \infty} \bar{E}\Big(\mathbb{P}(\Theta > N^{\gamma}t)\Big) = \lim\limits_{N \rightarrow \infty} \Bigg(1 - \frac{1}{N^{\gamma}}\Bigg)^{N^{\gamma}t+1} = e^{-t}.
	\]	 	
\end{proof}

\begin{proof}[Proof of Theorem \ref{T4}]
	Applying Chebyshev's inequality,
	
	\begin{align*}
	\bar{\mathbb{P}}\Big(|\mathbb{P}(\Theta > N^{\gamma}t) -  e^{-t}| > \epsilon\Big) \leq \! \frac{1}{\epsilon^{2}} \! \Bigg[ \! \bar{E} \!\Big( \! \Big[ \! \mathbb{P}(\Theta \! > \! N^{\gamma}t)\Big]^{2} \! \Big) \! - \! 2e^{-t} \bar{E}\Big(\mathbb{P}(\theta > N^{\gamma}t) \Big)+ e^{-2t}\Bigg].
	\end{align*}
	
	For a given $\bar{\omega} \in \bar{\Omega}$,
	
	\begin{align*}
	\Big(\mathbb{P}(\Theta > N^{\gamma}t)\Big)^{2} = & \Bigg\{\mathds{1}_{\{\eta \notin M_{N}\}} \frac{1}{N} \sum\limits_{i_{1} = 1}^{N} \Bigg(\mathds{1}_{\{\eta^{i_{1}} \notin M_{N}\}} \dots \frac{1}{N} \sum\limits_{i_{N^{\gamma}t} = 1}^{N} \mathds{1}_{\{\eta^{i_{1}, \dots, i_{N^{\gamma}t}} \notin M_{N}\}} \Bigg) \Bigg\} \times \\
	\times & \Bigg\{\mathds{1}_{\{\eta \notin M_{N}\}} \frac{1}{N} \sum\limits_{i^{*}_{1} = 1}^{N} \Bigg(\mathds{1}_{\{\eta^{i^{*}_{1}} \notin M_{N}\}} \dots \frac{1}{N} \sum\limits_{i^{*}_{N^{\gamma}t} = 1}^{N} \mathds{1}_{\{\eta^{i^{*}_{1}, \dots, i^{*}_{N^{\gamma}t}} \notin M_{N}\}} \Bigg) \Bigg\}.
	\end{align*}
	
	Let $G = \{\eta^{i^{*}_{1}}, \dots, \eta^{i^{*}_{1}, \dots, i^{*}_{N^{\gamma}t}}\}$. By Corollary \ref{col1},
	\[
	\lim\limits_{N \rightarrow \infty} |G^{*}| - N^{\gamma}t = 0.
	\]
	
	On the other hand, for a given $\{i^{*}_{1}, \dots, i^{*}_{N^{\gamma}t}\}$, by Proposition \ref{prop1},	
	\[
	\lim\limits_{N \rightarrow \infty} \mathbb{P}\Big(\{\eta,\eta^{i_{1}}, \dots, \eta^{i_{1}, \dots, i_{N^{\gamma}t}}\} \cap \{\eta^{i^{*}_{1}}, \dots, \eta^{i^{*}_{1}, \dots, i^{*}_{N^{\gamma}t}}\} \neq \emptyset\Big) = 0.
	\]
	
	Therefore,
	\[
	\bar{E}\Big[\Big(\mathbb{P}(\Theta > N^{\gamma}t)\Big)^{2}\Big] = \Bigg(1 - \frac{1}{N^{\gamma}}\Bigg)\Bigg(1 - \frac{1}{N^{\gamma}}\Bigg)^{2|G^{*}|} + o(N).
	\]
	
	When N diverges,	
	\[
	\lim\limits_{N \rightarrow \infty} \bar{E}\Big[\Big(\mathbb{P}(\Theta > N^{\gamma}t)\Big)^{2}\Big] = e^{-2t}.
	\]
	
	Thus, by Proposition \ref{prop5},
	\[
	\bar{\mathbb{P}}\Big(|\mathbb{P}(\Theta > N^{\gamma}t) -  e^{-t}| > \epsilon\Big) = 0.
	\]		
\end{proof}

\section*{Acknowledgements}
We would like to thank Antonio Galves for his orientation on the master's thesis \cite{tese}, in which this paper is based.

\bibliographystyle{alea3}
\bibliography{REF}

\begin{thebibliography}{14}
\providecommand{\natexlab}[1]{#1}
\providecommand{\url}[1]{\texttt{#1}}
\providecommand{\urlprefix}{URL }
\expandafter\ifx\csname urlstyle\endcsname\relax
  \providecommand{\doi}[1]{doi:\discretionary{}{}{}#1}\else
  \providecommand{\doi}{doi:\discretionary{}{}{}\begingroup
  \urlstyle{rm}\Url}\fi
\providecommand{\eprint}[2][]{\url{#2}}

\bibitem[{Cassandro et~al.(1984)Cassandro, Galves, Olivieri and Vares}]{tese2}
Marzio Cassandro, Antonio Galves, Enzo Olivieri and Maria~Eul{\'a}lia Vares.
\newblock Metastable behavior of stochastic dynamics: a pathwise approach.
\newblock \emph{Journal of statistical physics} \textbf{35}~(5-6), 603--634
  (1984).

\bibitem[{Cassandro et~al.(1991)Cassandro, Galves and Picco}]{tese1}
Marzio Cassandro, Antonio Galves and Pierre Picco.
\newblock Dynamical phase transitions in disordered systems: the study of a
  random walk model.
\newblock In \emph{Annales de l'IHP Physique th{\'e}orique}, volume~55, pages
  689--705 (1991).

\bibitem[{Cooper and Frieze(2014)}]{cooper}
Colin Cooper and Alan Frieze.
\newblock A note on the vacant set of random walks on the hypercube and other
  regular graphs of high degree.
\newblock \emph{Moscow Journal of Combinatorics and Number Theory}
  \textbf{4}~(4), 403--426 (2014).

\bibitem[{Crowe(1956)}]{crowe}
DW~Crowe.
\newblock The n-dimensional cube and the tower of hanoi.
\newblock \emph{The American Mathematical Monthly} \textbf{63}~(1), 29--30
  (1956).

\bibitem[{Diaconis et~al.(1990)Diaconis, Graham and Morrison}]{asy}
Persi Diaconis, Ronald~L. Graham and John~A. Morrison.
\newblock Asymptotic analysis of a random walk on a hypercube with many
  dimensions.
\newblock \emph{Random structures and algorithms} \textbf{1}~(1), 51--72
  (1990).

\bibitem[{Gilbert(1958)}]{gilbert}
Edgard~N Gilbert.
\newblock Gray codes and paths on the n-cube.
\newblock \emph{Bell Labs Technical Journal} \textbf{37}~(3), 815--826 (1958).

\bibitem[{Letac and Takacs(1979)}]{letac}
G.~Letac and L.~Takacs.
\newblock Random walks on an $m$-dimensional cube (stma v22 1919).
\newblock \emph{Journal für die Reine und Angewandte Mathematik} \textbf{310},
  187--195 (1979).

\bibitem[{Matthews(1987)}]{mat2}
Peter Matthews.
\newblock Mixing rates for a random walk on the cube.
\newblock \emph{SIAM Journal on Algebraic Discrete Methods} \textbf{8}~(4),
  746--752 (1987).

\bibitem[{Matthews(1989)}]{mat}
Peter Matthews.
\newblock Some sample path properties of a random walk on the cube.
\newblock \emph{Journal of Theoretical Probability} \textbf{2}~(1), 129--146
  (1989).

\bibitem[{Nestoridi(2017)}]{nestoridi2017}
Evita Nestoridi.
\newblock A non-local random walk on the hypercube.
\newblock \emph{Advances in Applied Probability} \textbf{49}~(4), 1288--1299
  (2017).

\bibitem[{Peixoto(1992)}]{tese}
Cl\'adia Peixoto.
\newblock \emph{{Aproxima\c{c}\~ao do Equil\'ibrio e Tempos Exponenciais para o
  Passeio Aleat\'orio no Hipercubo}}.
\newblock Master's thesis, Instituto de Matem\'atica e Estat\'istica da
  Universidade de S\~ao Paulo, S\~ao Paulo, Brasil. (1992).
\newblock Available at \url{https://arxiv.org/abs/1805.11349}.

\bibitem[{Scoppola(2011)}]{sco}
Benedetto Scoppola.
\newblock Exact solution for a class of random walk on the hypercube.
\newblock \emph{Journal of Statistical Physics} \textbf{143}~(3), 413--419
  (2011).

\bibitem[{Voit(1996)}]{ehr}
Michael Voit.
\newblock Asymptotic distributions for the ehrenfest urn and related random
  walks.
\newblock \emph{Journal of Applied Probability} \textbf{33}~(2), 340--356
  (1996).
\newblock ISSN 00219002.
\newblock \urlprefix\url{http://www.jstor.org/stable/3215058}.

\bibitem[{Volkov and Wong(2008)}]{volkov}
Stanislav Volkov and Timothy Wong.
\newblock A note on random walks in a hypercube.
\newblock \emph{Pi Mu Epsilon journal} pages 551--557 (2008).

\end{thebibliography}

\end{document}